\documentclass[a4paper,10pt]{article}
\usepackage[utf8x]{inputenc}
\usepackage{tracefnt,amsmath,tabu,array}
\usepackage{amssymb,graphicx,setspace,amsfonts,amsbsy}
\usepackage{pifont,latexsym,ifthen,amsthm,rotating,calc,textcase,booktabs,cancel,slashed}

\newtheorem{theorem}{Theorem}[section]
\newtheorem{lemma}[theorem]{Lemma}
\newtheorem{corollary}[theorem]{Corollary}
\newtheorem{proposition}[theorem]{Proposition}
\newtheorem{definition}[theorem]{Definition}

\newtheorem{remark}[theorem]{Remark}
\newcommand{\filledbox}{\leavevmode
  \hbox to.77778em{%
  \hfil\vbox to.675em{\hrule width.6em height.6em}\hfil}}

\newcommand{\Rm}{{\mathbb R}}

\begin{document}
\tabulinesep=1.0mm
\title{Radiation fields and non-radiative solutions to the energy sub-critical wave equations}

\author{Liang Li, Ruipeng Shen and Chenhui Wang\\
Centre for Applied Mathematics\\
Tianjin University\\
Tianjin, China
}

\maketitle

\begin{abstract}
 Radiation field and channel of energy method have become important tools in the study of nonlinear wave equations in recent years. In this work we give basic theory of radiation fields of free waves in the energy sub-critical case. We also show that the asymptotic behaviours of non-radiative solutions to a wide range of non-linear wave equations resemble those of non-radiative free waves. Our theory is completely given in the critical Sobolev spaces of the corresponding nonlinear wave equation and avoids any assumption on the energy of the solutions. 
\end{abstract}

\section{Introduction}

In this work we discuss the properties of radiation fields of wave equations in energy sub-critical spaces and generalize the channel of energy method to these spaces. 

\subsection{Background} 

\paragraph{Radiation fields} The history of radiation fields is more than 50 years long. Please see, Friedlander \cite{radiation1, radiation2}, for instance. In their earlier history radiation fields were mainly conceptions of mathematical physics. In recent years the radiation fields play an important role in the study of asymptotic behaviour of solutions to non-linear wave equations. Briefly speaking, radiation fields describe the asymptotic behaviour of linear free waves, i.e. the solution to the homogenous linear wave equation 
\[
 \partial_t^2 u - \Delta u = 0.
\]
The following version of radiation fields is given by Duyckaerts-Kenig-Merle \cite{dkm3}.

\begin{theorem}[Radiation field] \label{radiation}
Assume that $d\geq 3$ and let $u$ be a solution to the free wave equation $\partial_t^2 u - \Delta u = 0$ with initial data $(u_0,u_1) \in \dot{H}^1 \times L^2(\Rm^d)$. Then ($u_r$ is the derivative in the radial direction)
\[
 \lim_{t\rightarrow \pm \infty} \int_{\Rm^d} \left(|\nabla u(x,t)|^2 - |u_r(x,t)|^2 + \frac{|u(x,t)|^2}{|x|^2}\right) dx = 0
\]
 and there exist two functions $G_\pm \in L^2(\Rm \times \mathbb{S}^{d-1})$ so that
\begin{align*}
 \lim_{t\rightarrow \pm\infty} \int_0^\infty \int_{\mathbb{S}^{d-1}} \left|r^{\frac{d-1}{2}} \partial_t u(r\theta, t) - G_\pm (r\mp t, \theta)\right|^2 d\theta dr &= 0;\\
 \lim_{t\rightarrow \pm\infty} \int_0^\infty \int_{\mathbb{S}^{d-1}} \left|r^{\frac{d-1}{2}} \partial_r u(r\theta, t) \pm G_\pm (r\mp t, \theta)\right|^2 d\theta dr & = 0.
\end{align*}
In addition, the maps $(u_0,u_1) \rightarrow \sqrt{2} G_\pm$ are bijective isometries from $\dot{H}^1 \times L^2(\Rm^d)$ to $L^2 (\Rm \times \mathbb{S}^{d-1})$. 
\end{theorem}

\paragraph{Radiation profiles} We call the function $G_\pm$ given above radiation profiles of the initial data $(u_0,u_1)$ (or equivalently, of the linear free wave $u$). We also define $\mathbf{T}_\pm$ to be the map $(u_0,u_1)\rightarrow G_\pm$. These maps can be given in term of Radon transform (see Katayama \cite{katayamaradiation}).

\paragraph{Radiation fields by Fourier transforms} The maps $\mathbf{T}_\pm$ can also be given explicitly in the form of Fourier transforms. We use the notation $\hat{u}$ for the Fourier transform in $\Rm^d$
\[
 \hat{u} (\xi) = \int_{\Rm^d} e^{i\xi\cdot x} u(x) dx. 
\]  
and $\mathcal{F} G$ for the partial Fourier transform in the first variable $s \in \Rm$ 
\[
 (\mathcal{F} G) (\nu, \omega) = \int_{\Rm} e^{i\nu s} G(\nu, \omega) ds. 
\]
The radiation fields $G_\pm = \mathbf{T}_\pm (u_0,u_1)$ can be given in term of the initial data (see C\"{o}te-Laurent \cite{newradiation}, please note that our notations are slightly different from their original work)
\begin{align}
 \mathcal{F} G_+  & =   - c_d |\nu|^\frac{d-1}{2} \left[e^{i\tau} \chi_-(\nu)(i\nu \hat{u}_0 (\nu \omega)-\hat{u}_1(\nu \omega))+e^{-i\tau}\chi_+(\nu)(i\nu \hat{u}_0 (\nu \omega)-\hat{u}_1(\nu \omega))\right]; \label{formula G plus} \\
 \mathcal{F} G_-  & = + c_d  |\nu|^\frac{d-1}{2}\left[e^{i\tau} \chi_-(\nu)(i\nu \hat{u}_0 (\nu \omega)+\hat{u}_1(\nu \omega))+e^{-i\tau}\chi_+(\nu) (i\nu \hat{u}_0 (\nu \omega)+\hat{u}_1(\nu \omega))\right]. \label{formula G negative}
\end{align}
Here 
\begin{align*}
 &c_d = \frac{1}{2(2\pi)^{\frac{d-1}{2}}},& &\tau = \frac{d-1}{4}\pi.&
\end{align*}
Their inverse $\mathbf{T}_\pm^{-1}$ can also given explicitly in the form of Fourier transforms by a basic calculation.
\begin{align*}
 \hat{u}_0 (\nu \omega) & = \frac{i}{2c_d} |\nu|^{-\frac{d+1}{2}}\left[e^{i\tau} (\mathcal{F} G_+) (\nu, \omega) - e^{-i\tau} (\mathcal{F} G_+)(-\nu, -\omega)\right];\\
 \hat{u}_1 (\nu \omega) & = \frac{1}{2c_d} |\nu|^{-\frac{d-1}{2}}\left[e^{i\tau} (\mathcal{F} G_+) (\nu, \omega) + e^{-i\tau} (\mathcal{F} G_+)(-\nu, -\omega)\right];\\
  \hat{u}_0 (\nu \omega) & = \frac{-i}{2c_d} |\nu|^{-\frac{d+1}{2}}\left[e^{i\tau} (\mathcal{F} G_-) (\nu, \omega) - e^{-i\tau} (\mathcal{F} G_-)(-\nu, -\omega)\right];\\
 \hat{u}_1 (\nu \omega) & = \frac{1}{2c_d} |\nu|^{-\frac{d-1}{2}}\left[e^{i\tau} (\mathcal{F} G_-) (\nu, \omega) + e^{-i\tau} (\mathcal{F} G_-)(-\nu, -\omega)\right].
\end{align*}

\paragraph{Explicit formula of $\mathbf{T}_\pm^{-1}$} Besides the expression in the form of Fourier transform, we may also give an explicit formula of the inverse $\mathbf{T}_\pm^{-1}: G_\pm \rightarrow (u_0,u_1)$ in term of integral transformation. The following formula comes from Li-Shen-Wei \cite{shenradiation}.  If $d \geq 3$ is odd, then we have
\begin{align*}
  u_0(x) & = \frac{1}{(2\pi)^\frac{d-1}{2}} \int_{\mathbb{S}^{d-1}}  \partial_s^\frac{d-3}{2} G_- \left(x \cdot \omega, \omega\right) d\omega\\
 u_1(x) & = \frac{1}{(2\pi)^\frac{d-1}{2}} \int_{\mathbb{S}^{d-1}}  \partial_s^\frac{d-1}{2} G_- \left(x \cdot \omega, \omega\right)  d\omega
\end{align*}
If $d\geq 2$ is even, then
\begin{align*}
 u_0 (x) & = \frac{\sqrt{2}}{(2\pi)^{d/2}} \cdot  \int_0^{\infty} \int_{\mathbb{S}^{d-1}} \frac{\partial_s^{d/2-1} G_- \left(x\cdot \omega -\rho, \omega\right)}{\sqrt{\rho}} d\omega d\rho;\\
 u_1 (x) & =\frac{\sqrt{2}}{(2\pi)^{d/2}} \cdot  \int_0^{\infty} \int_{\mathbb{S}^{d-1}} \frac{\partial_s^{d/2} G_- \left(x\cdot \omega -\rho, \omega\right)}{\sqrt{\rho}} d\omega d\rho.
\end{align*}
By a time translation, we immediately have
\begin{align*}
 u(x,t) & = \frac{1}{(2\pi)^\frac{d-1}{2}} \int_{\mathbb{S}^{d-1}}  \partial_s^\frac{d-3}{2} G_- \left(x \cdot \omega +t , \omega\right) d\omega, & &\hbox{$d$ is odd;}\\
 u(x,t) & = \frac{\sqrt{2}}{(2\pi)^{d/2}} \cdot  \int_0^{\infty} \int_{\mathbb{S}^{d-1}} \frac{\partial_s^{d/2-1} G_- \left(x\cdot \omega -\rho + t, \omega\right)}{\sqrt{\rho}} d\omega d\rho,& &\hbox{$d$ is even.} 
\end{align*}

\paragraph{Non-radiative solutions} We call a solution $u$ to wave equation ($R$-weakly) non-radiative if 
\[
 \lim_{t\rightarrow \pm \infty} \int_{|x|>|t|+R} |\nabla_{t,x} u(x,t)|^2 dx = 0.
\]
This kind of solutions are crucial in the channel of energy method, which becomes a powerful tool in the study of wave equations in the past decade. More details about the channel of energy method can be found in, or instance,  C\^{o}te-Kenig-Schlag \cite{channeleven}, Duyckaerts-Kenig-Merle \cite{tkm1, oddtool} and Kenig-Lawrie-Schlag \cite{channel5d}. Let $u$ be a free wave with a finite energy and $G_\pm$ be its radiation profiles. Then Theorem \ref{radiation} implies that $u$ is $R$-weakly non-radiative if and only if $G_\pm (s,\omega) = 0$ if $s>R$. If the space dimension $d$ is odd, then we have
\[
 G_+(s,\omega) = (-1)^{\frac{d-1}{2}} G_-(-s,-\omega)
\]
In fact, this is a direct consequence of the explicit formula \eqref{formula G plus} and \eqref{formula G negative}. Thus if $d$ is odd, a solution is $R$-weakly non-radiative solution if and only if $G_-$ is supported in the set $[-R,R]\times \mathbb{S}^{d-1}$. 

\paragraph{Non-linear equations} We then consider energy-critical non-linear wave equation in 3-dimensional space
\[
 \partial_t^2 u - \Delta u = F(x,t,u).
\]
Here the non-linear term $F$ satisfies 
\begin{align*}
 &|F(x,t,u)| \lesssim |u|^5;& &|F(x,t,u_1)-F(x,t,u_2)|\lesssim (|u_1|^4+|u_2|^4)|u_1-u_2|.&
\end{align*} 
By applying a decay estimate of linear non-radiative solutions in the exterior region $\{(x,t): |x|>|t|+r\}$ given in Li-Shen-Wang \cite{3Dnonradiativedecay}, we may obtain a decay estimate of the radiation profile $G_-$ of initial data for any given non-radiative solution $u$ in Li-Shen-Wang-Wei \cite{nonradialCE}: 
\[
 \|G_-(s,\omega)\|_{L^2(((-\infty,r]\cap [r,+\infty))\times \mathbb{S}^2)} \leq r^{-\kappa}, \qquad r\geq R(\kappa, u).
\]
Here $\kappa$ is any positive constant less than $1$. In other words, the radiation profiles associated to initial data of non-linear non-radiative solutions asymptoticly look like the radiation profiles of linear non-radiative solutions, which are compactly-supported. This result deals non-radial solutions as well, while most similar results are concerning radial solutions, see, Duyckaerts-Kenig-Merle\cite{se, oddtool}, for example.

\subsection{Main topics}

In this work we try to generalize the basic theory of radiation fields and related results about non-radiative solutions to the energy sub-critical space. Our results are motivated by the following basic observation
\begin{proposition} \label{prop isometry}
 Given $\beta \in \Rm$, the maps $\sqrt{2} \mathbf{T}_\pm$ are isometries from the space $\dot{H}^{\beta} (\Rm^d) \times \dot{H}^{\beta -1}(\Rm^d)$ to the space $L^2(\mathbb{S}^{d-1}, \dot{H}^{\beta-1}(\Rm))$. 
\end{proposition}
\begin{proof}
 We use the same notations of radiation profiles and Fourier transforms as given in the last subsection. Both $\mathbf{T}_\pm$ and their inverses have been given explicitly. It suffices to verify that $\sqrt{2} \mathbf{T}_\pm$ preserve the norm. A straight forward calculation shows 
\begin{align*}
 |(\mathcal{F} G_+)(\nu, \omega)|^2 & = c_d^2 |\nu|^{d-1} \left|i\nu \hat{u}_0 (\nu \omega)-\hat{u}_1(\nu \omega)\right|^2 \\
 & = c_d^2 |\nu|^{d-1} \left(\nu^2 |\hat{u}_0(\nu \omega)|^2 + |\hat{v}_1(\nu \omega)|^2 -2 \hbox{Re}(i \nu \hat{v}_0(\nu \omega) \hat{v}_1(\nu \omega)) \right).
\end{align*}
Thus
\[
 |(\mathcal{F} G_+)(-\nu, -\omega)|^2 = c_d^2 |\nu|^{d-1} \left(\nu^2 |\hat{u}_0(\nu \omega)|^2 + |\hat{v}_1(\nu \omega)|^2 +2 \hbox{Re}(i \nu \hat{v}_0(\nu \omega) \hat{v}_1(\nu \omega)) \right).
\]
Therefore 
\begin{align*}
 \sum_{\pm} \int_{\mathbb{S}^{d-1}} |(\mathcal{F} G_+)(\pm \nu, \omega)|^2 d\omega & = \int_{\mathbb{S}^{d-1}} \left(|(\mathcal{F} G_+)(\nu, \omega)|^2 + |(\mathcal{F} G_+)(-\nu, -\omega)|^2\right) d\omega \\
 & = 2c_d^2 \int_{\mathbb{S}^{d-1}} |\nu|^{d-1}\left(\nu^2 |\hat{u}_0(\nu \omega)|^2 + |\hat{v}_1(\nu \omega)|^2\right) d\omega
\end{align*}
This immediately gives ($\beta \in \Rm$)
\begin{equation} \label{isometry identity}
 2 \int_{\Rm\times \mathbb{S}^{d-1}} \left|(\mathbf{D}_s^{\beta-1} G_+)(s,\omega)\right|^2 ds d\omega = \|(u_0,u_1)\|_{\dot{H}^{\beta}(\Rm^d) \times \dot{H}^{\beta-1}(\Rm^d)}^2
\end{equation}
Here $\mathbf{D}_s^\alpha$ is the (fractional) derivative with respect to $s$, defined by
\[
 \mathcal{F} (\mathbf{D}_s^\alpha G) (\nu, \omega) = |\nu|^\alpha \mathcal{F} G (\nu, \omega).
\]
This finishes the proof in the positive time direction. The negative direction is similar. 
\end{proof}

\paragraph{The goal} In this work we will give two main results, both in the energy sub-critical setting. 
\begin{itemize}
 \item We consider linear free waves $u$ with initial data $(u_0,u_1) \in \dot{H}^{\beta} \times \dot{H}^{\beta-1}(\Rm^d)$. Here $\beta \in (0,1)$ and $d\geq 3$. Let $G_\pm$ be the radiation profiles of $u$. Then we have 
  \begin{align*}
  \lim_{t\rightarrow +\infty} \left\|\nabla_{t,x} u(x,t) - \frac{1}{|x|^{\frac{d-1}{2}}} G_+ \left(|x|-t,\frac{x}{|x|}\right)\times \begin{pmatrix} 1\\ -x/|x|\end{pmatrix}\right\|_{(\dot{H}^{\beta-1}(\Rm^d))^{1+d}} & = 0;\\
  \lim_{t\rightarrow -\infty} \left\|\nabla_{t,x} u(x,t) - \frac{1}{|x|^{\frac{d-1}{2}}} G_- \left(|x|+t,\frac{x}{|x|}\right)\times \begin{pmatrix} 1\\ x/|x|\end{pmatrix}\right\|_{(\dot{H}^{\beta-1}(\Rm^d))^{1+d}} & = 0.
 \end{align*}
 \item We consider a wide range of non-linear wave equation in dimension $3$ whose critical Sobolev space is $\dot{H}^{s_p} \times \dot{H}^{s_p-1} (\Rm^3)$ with $s_p \in (1/2,1)$. We show that if $u$ is $R$-weakly non-radiative, then the radiation profiles of its initial data must decay at a certain rate as $s\rightarrow \infty$. Thus the radiation profiles asymptoticly look like those of linear non-radiative solutions. Please note that the $R$-weakly non-radiative solutions here are defined in a way that we may completely work in the space $\dot{H}^{s_p} \times \dot{H}^{s_p-1} (\Rm^d)$, thus get rid of the use of energy. More details can be found in Section \ref{sec: non-linear wave equation}. 
\end{itemize}

\paragraph{Structure of this work} In Section 2 we prove the first main result given above, i.e. radiation fields of linear free waves in the energy sub-critical space. Next we give a few preliminary results and decay estimates regarding non-radiative linear free waves in the energy sub-critical space in Section 3. Finally in Section 4 we prove our second main result regarding non-radiative non-linear solutions. 
 
\paragraph{Notations} Before we conclude this section we give a few notations, which will be used throughout this paper. 
\begin{align*}
 \dot{H}^\beta (\Rm \times \mathbb{S}^{d-1}) &\doteq L^2(\mathbb{S}^{d-1}; \dot{H}^\beta (\Rm)); \\
 \dot{H}_0^\beta (I \times \mathbb{S}^{d-1}) & \doteq L^2(\mathbb{S}^{d-1}; \dot{H}_0^\beta (I)), \qquad \beta \in (-1/2,1/2);\\
 \mathcal{H}^{\beta} (\Rm^d) & \doteq \dot{H}^{\beta}(\Rm^d) \times \dot{H}^{\beta-1}(\Rm^d).  
\end{align*}
Here $I$ is an interval and $\dot{H}_0^\beta (I)$ is the completion of $C_0^\infty(I)$ in $\dot{H}^\beta (\Rm)$. Given initial data $(u_0,u_1)$, we also use the notation $\mathbf{S}_L (u_0,u_1)$ for the solution to the free linear wave equation with initial data $(u_0,u_1)$. In this work we also use the notation $A \lesssim B$ to indicate that there exists a constant $c$ so that $A \leq cB$. We may add subscripts to the notation $\lesssim$ to indicate that the explicit constant depends on these subscripts but nothing else. 

\section{Radiation fields in energy sub-critical space}
\begin{lemma} \label{dual lemma}
 Let $\beta \in [0,1]$ be a constant and $d \geq 3$. Then 
 \begin{itemize}
  \item[(i)] The following linear operators $\mathbf{A}_k^\ast (k=0,1,\cdots, d)$ from $\dot{H}^\beta (\Rm^d)$ to $\dot{H}^\beta (\Rm \times \mathbb{S}^{d-1})$ is uniformly bounded for all $t > 0$.
  \[
   f(x) \rightarrow (\mathbf{A}_k^\ast f)(s,\omega) = \left\{\begin{array}{ll} (s+t)^\frac{d-1}{2} f((s+t)\omega) \omega_k, & s>-t; \\ 0, & s\leq -t. \end{array} \right.
  \]
  Here $\omega = (\omega_1, \omega_2, \cdots, \omega_d) \in \mathbb{S}^{d-1}$ and $\omega_0 = 1$.
  \item[(ii)] The linear operator $\mathbf{A}$ from $\dot{H}^{-\beta} (\Rm \times \mathbb{S}^{d-1})$ to $(\dot{H}^{-\beta} (\Rm^d))^{1+d}$ is uniformly bounded for all $t>0$.
  \[
   G(s,\omega) \rightarrow (\mathbf{A} G) (x) = \frac{1}{|x|^\frac{d-1}{2}} G\left(|x|-t, \frac{x}{|x|}\right) \times \left(\begin{array}{c} 1 \\ -x/|x|\end{array}\right).
  \]
 \end{itemize}
\end{lemma}
\begin{proof}
 First of all, we may utilize polar coordinates to verify
 \begin{align*}
  \int_{\Rm^d} f(x) \cdot \frac{1}{|x|^\frac{d-1}{2}} G\left(|x|-t, \frac{x}{|x|}\right) dx &= \int_{(-t,+\infty) \times \mathbb{S}^{d-1}} (s+t)^\frac{d-1}{2} f((s+t)\omega)  G(s, \omega) ds d\omega; \\
  \int_{\Rm^d} f(x) \cdot \frac{x}{|x|^\frac{d+1}{2}} G\left(|x|-t, \frac{x}{|x|}\right) dx &= \int_{(-t,+\infty) \times \mathbb{S}^{d-1}} (s+t)^\frac{d-1}{2} f((s+t)\omega)  G(s, \omega) \omega ds d\omega.
 \end{align*}
 Thus it suffices to prove part (i) by duality. First of all, a basic calculation shows that $\mathbf{A}_k^\ast$ is a bounded from $L^2(\Rm^d)$ to $L^2(\Rm \times \mathbb{S}^{d-1})$. This deals with the case $\beta = 0$. Next we consider the case $\beta = 1$. We may calculate 
 \begin{align*}
  \int_{(-t,+\infty) \times \mathbb{S}^{d-1}} \left|\frac{d}{ds} \left[(s+t)^\frac{d-1}{2} f((s+t)\omega) \omega_k \right]\right|^2 ds d\omega & \leq \int_{\mathbb{S}^{d-1}} \int_0^\infty \left|\frac{(d-1)f}{2r}+f_r\right|^2 r^{d-1} dr d\omega \\
  & \lesssim_d \|f\|_{\dot{H}^1(\Rm^d)}^2. 
 \end{align*}
Here we apply the Hardy's inequality. This gives the uniform boundedness of $\mathbf{A}_k^\ast$ in the case $\beta = 1$. Finally we apply an interpolation between $\beta = 0$ and $\beta = 1$ to finish our proof. 
\end{proof}
\begin{proposition}[Radiation fields in energy sub-critical space]
 Assume that $(u_0,u_1)\in \dot{H}^{\beta} \times \dot{H}^{\beta-1}(\Rm^d)$ with $\beta \in (0,1)$. Let $G_\pm = \mathbf{T}_\pm (u_0,u_1)$ be the radiation profiles of $(u_0,u_1)$. Then the free wave $u =\mathbf{S}_L (u_0,u_1)$ satisfies 
 \begin{align*}
  \lim_{t\rightarrow +\infty} \left\|\nabla_{t,x} u(x,t) - \frac{1}{|x|^{\frac{d-1}{2}}} G_+ \left(|x|-t,\frac{x}{|x|}\right)\times \begin{pmatrix} 1\\ -x/|x|\end{pmatrix}\right\|_{(\dot{H}^{\beta-1}(\Rm^d))^{1+d}} & = 0;\\
  \lim_{t\rightarrow -\infty} \left\|\nabla_{t,x} u(x,t) - \frac{1}{|x|^{\frac{d-1}{2}}} G_- \left(|x|+t,\frac{x}{|x|}\right)\times \begin{pmatrix} 1\\ x/|x|\end{pmatrix}\right\|_{(\dot{H}^{\beta-1}(\Rm^d))^{1+d}} & = 0.
 \end{align*}
\end{proposition}
\begin{proof}
 Let us prove the limit in the positive time direction. The other limit can be proved in the same way. A combination of Proposition \ref{prop isometry} and Lemma \ref{dual lemma}, as well as the fact that the linear wave propagation operator preserves $\dot{H}^\beta \times \dot{H}^{\beta-1}$ norm, shows that  
 \[
  \left\|\nabla_{t,x} u(x,t) - \frac{1}{|x|^{\frac{d-1}{2}}} G_+ \left(|x|-t,\frac{x}{|x|}\right)\times \begin{pmatrix} 1\\ -x/|x|\end{pmatrix}\right\|_{(\dot{H}^{\beta-1}(\Rm^d))^{1+d}} \lesssim_d \|(u_0,u_1)\|_{\dot{H}^\beta \times \dot{H}^{\beta-1}}. 
 \]
 Thus it suffices to verity the limit for initial data in a dense subset of $\dot{H}^{\beta} \times \dot{H}^{\beta-1}$ by the linearity. We fix $\beta_0 \in (0,\beta)$ and consider initial data in the space $(\dot{H}^1 \cap \dot{H}^{\beta_0}) \times (L^2 \cap \dot{H}^{\beta_0-1})$. In this case the radiation field in the energy space implies 
 \[
  \lim_{t\rightarrow +\infty} \left\|\nabla_{t,x} u(x,t) - \frac{1}{|x|^{\frac{d-1}{2}}} G_+ \left(|x|-t,\frac{x}{|x|}\right)\times \begin{pmatrix} 1\\ -x/|x|\end{pmatrix}\right\|_{(L^2(\Rm^d))^{1+d}} = 0.
 \]
 In addition we have
 \[
  \left\|\nabla_{t,x} u(x,t) - \frac{1}{|x|^{\frac{d-1}{2}}} G_+ \left(|x|-t,\frac{x}{|x|}\right)\times \begin{pmatrix} 1\\ -x/|x|\end{pmatrix}\right\|_{(\dot{H}^{\beta_0-1}(\Rm^d))^{1+d}} \lesssim_d \|(u_0,u_1)\|_{\dot{H}^{\beta_0} \times \dot{H}^{\beta_0-1}}. 
 \]
 Finally we may finish the proof by an interpolation. 
\end{proof}

\section{Decay estimates of linear waves}
  
\subsection{Notations and preliminary results}

We first introduce a few basic notations. 

\paragraph{Cut-off operators} Let $\chi_1$ and $\chi_2$ be the characteristic functions of the sets $(1,+\infty)$ and $(0,1)$, respectively. 
Then both the maps $f(s) \rightarrow \chi_1 (s) f(s)$ and $f(s) \rightarrow \chi_2 (s) f(s)$ are bounded operators from $\dot{H}^\beta(\Rm)$ to itself if $\beta \in (-1/2,1/2)$. Please see Triebel \cite{cutoff}, for instance. We may apply translations, dilations and reflections to conclude that the cut-off operators defined below are all bounded operators from $\dot{H}^\beta(\Rm)$ to itself
\begin{align*} 
 &\mathbf{P}_r^\pm f = \chi_1 (\pm s/r) f(s);& &\mathbf{P}_{a,b} f = \chi_2\left(\frac{s-a}{b-a}\right) f(s).&
\end{align*}
In addition, the operator norms of these operators are independent to the constants $r > 0$, $a, b \in \Rm$. We may also rewrite these operators in the following way
\begin{align*}
 &(\mathbf{P}_r^+ f)(s) = \left\{\begin{array}{ll} f(s), & s>r; \\ 0, & s\leq r; \end{array}\right.& &(\mathbf{P}_r^- f)(s) = \left\{\begin{array}{ll} f(s), & s<-r; \\ 0, & s\geq -r; \end{array}\right.& \\
 &(\mathbf{P}_{a,b} f)(s) = \left\{\begin{array}{ll} f(s), & a<s<b; \\ 0, & \hbox{otherwise}. \end{array}\right.& &&
\end{align*}
These operators can be applied on $G \in \dot{H}^\beta (\Rm\times \mathbb{S}^2)$ in a naturally way 
\begin{align*} 
 &(\mathbf{P}_r^\pm f)(s,\omega) = \chi_1 (\pm s/r) f(s, \omega);& &(\mathbf{P}_{a,b} f)(s, \omega) = \chi_2 \left(\frac{s-a}{b-a}\right) f(s, \omega).&
\end{align*}
The operators $\mathbf{P}_r^+$, $\mathbf{P}_r^-$ and $\mathbf{P}_{a,b}$ are bounded operators from $\dot{H}^\beta (\Rm \times \mathbb{S}^{d-1})$ to $\dot{H}_0^\beta ((r,+\infty)\times \mathbb{S}^{d-1})$, $\dot{H}_0^\beta ((-\infty,-r)\times \mathbb{S}^{d-1})$ and $\dot{H}_0^\beta ((a,b)\times \mathbb{S}^{d-1})$, respectively. Here we assume $\beta \in (-1/2, 1/2)$. 

\paragraph{Strichartz estimates} We next recall the generalized Strichartz estimates. Please see Proposition 3.1 in Ginibre-Velo \cite{strichartz}. Here we use the Sobolev version in dimension 3.
\begin{proposition} [Generalized Strichartz estimates] Let $2\leq q_1,q_2 \leq \infty$, $2\leq r_1,r_2 < \infty$ and $s\in \Rm$ be constants with 
\begin{align*}
 &1/q_i+1/r_i \leq 1/2, \; i=1,2; & &1/q_1+3/r_1=3/2-s;& &1/q_2+3/r_2=1/2+s.&
\end{align*}
Assume that $u$ is the solution to the linear wave equation
\[
 \left\{\begin{array}{ll} \partial_t u - \Delta u = F(x,t), & (x,t) \in \Rm^3 \times [0,T];\\
 u|_{t=0} = u_0 \in \dot{H}^s; & \\
 \partial_t u|_{t=0} = u_1 \in \dot{H}^{s-1}. &
 \end{array}\right.
\]
Then we have
\begin{align*}
 \left\|\left(u(\cdot,T), \partial_t u(\cdot,T)\right)\right\|_{\dot{H}^s \times \dot{H}^{s-1}} & +\|u\|_{L^{q_1} L^{r_1}([0,T]\times \Rm^3)} \\
 & \leq C\left(\left\|(u_0,u_1)\right\|_{\dot{H}^s \times \dot{H}^{s-1}} + \left\|F(x,t) \right\|_{L^{\bar{q}_2} L^{\bar{r}_2} ([0,T]\times \Rm^3)}\right).
\end{align*}
Here the coefficients $\bar{q}_2$ and $\bar{r}_2$ satisfy $1/q_2 + 1/\bar{q}_2 = 1$, $1/r_2 + 1/\bar{r}_2 = 1$. The constant $C$ does not depend on $T$ or $u$. 
\end{proposition}
\noindent In particular we introduce the following notations 
\begin{align*}
 &\|u\|_{Y (I)} = \|u\|_{L^\frac{2p}{s_p+1} L^\frac{2p}{2-s_p}(I \times \Rm^3)}; & &\|u\|_{Z(I)} = \|u\|_{L^\frac{2}{s_p+1} L^\frac{2}{2-s_p}(I \times \Rm^3)}.&
\end{align*}
Here $s_p = 3/2 - 2/(p-1)$ is the critical exponent of the non-linear wave equation 
\[
 \partial_t^2 u - \Delta u = \pm |u|^{p-1} u.
\]
Then the Strichartz estimate given above becomes 
\begin{equation} \label{strichartz YZ}
 \left\|\left(u(\cdot,T), \partial_t u(\cdot,T)\right)\right\|_{\dot{H}^{s_p} \times \dot{H}^{s_p-1}} +\|u\|_{Y([0,T])}  \leq C\left(\left\|(u_0,u_1)\right\|_{\dot{H}^{s_p} \times \dot{H}^{s_p-1}} + \left\|F(x,t) \right\|_{Z ([0,T])}\right). 
\end{equation}
We conclude this subsection by a technical lemma 

\begin{lemma}
 Assume that $q>0$, $\gamma>\gamma_1 > 1$ and $a_j \in [0,+\infty)$. Then  
 \[
  \left(\sum_{j=0}^{+\infty} \gamma^{-j} a_j\right)^q \leq \left(\frac{\gamma}{\gamma-\gamma_1}\right)^q \left(\sum_{j=0}^{+\infty} \gamma_1^{-jq} a_j^q\right)
 \]
\end{lemma}
\begin{proof}
 This follows a basic calculation 
 \begin{align*}
  \sum_{j=0}^{+\infty} \gamma^{-j} a_j & =  \sum_{j=0}^{+\infty} (\gamma_1/\gamma)^j \gamma_1^{-j} a_j \leq \left(\sum_{j=0}^{+\infty} (\gamma_1/\gamma)^j \right)\left(\sum_{j=0}^{+\infty} \gamma_1^{-jq} a_j^q \right)^{1/q} \leq \frac{\gamma}{\gamma-\gamma_1} \left(\sum_{j=0}^{+\infty} \gamma_1^{-jq} a_j^q \right)^{1/q}.
 \end{align*}
\end{proof}

\subsection{Decay estimates}

In this subsection we consider the decay estimate regarding the operator 
\[
 \mathbf{T} G (x)= \int_{\mathbb{S}^2} G(x\cdot \omega, \omega) d\omega. 
\]
This is exactly the adjoint of Radon transform. We recall the following $L^6$ decay estimates when the radiation profile $G \in L^2 (\Rm \times \mathbb{S}^2)$ is compactly supported in $[a,b]\times \mathbb{S}^2$, as given in Li-Shen-Wang \cite{3Dnonradiativedecay}. 
\begin{equation} \label{left point interpolation}
  \int_{|x|>R} |\mathbf{T} G (x)|^6 dx \lesssim  \frac{(b-a)^2}{R^2}  \|G\|_{L^2(\Rm\times \mathbb{S}^2)}^6, \qquad R>0.
\end{equation} 
If $G$ is supported in $([-b,-a] \cup[a,b])\times \mathbb{S}^2$ with $b>a>0$, then we have $\mathbf{T} G (x) = 0$ if $|x|<a$, thus we also have the global $L^6$ estimate 
\begin{equation} \label{left point interpolation 2}
  \int_{\Rm^3} |\mathbf{T} G (x)|^6 dx \leq \int_{|x|>a} |\mathbf{T} G (x)|^6 dx \lesssim  (b/a-1)^2  \|G\|_{L^2(\Rm\times \mathbb{S}^2)}^6.
\end{equation} 
Now we prove 
\begin{lemma} \label{Y estimate}
 Assume that $p \in (3,5)$, $I = (a,b)$ and $G \in \dot{H}_0^{s_p-1}(I \times \mathbb{S}^2)$. Then the linear free wave (whose radiation profile is exactly $G$)
 \[
  u(x,t) = \frac{1}{2\pi} \int_{\mathbb{S}^2} G(x\cdot \omega+t, \omega) d\omega
 \]
 satisfies ($s_p = 3/2-2/(p-1)$)
 \[
  \|u\|_{L_t^\frac{2p}{1+s_p} L^\frac{2p}{2-s_p}(\{x: |x|>R\})} \lesssim \left(\frac{b-a}{R}\right)^{\frac{p-3}{2p}} \|G\|_{\dot{H}_0^{s_p-1}(I \times \mathbb{S}^2)}.
 \]
\end{lemma}
\begin{proof}
 Let us consider more general radiation profiles, i.e. those without the compact support assumption, and define a map by
 \[
  G \rightarrow u(x,t) \doteq \mathbf{S}_L \mathbf{T}_-^{-1} \mathbf{P}_{a,b} G = \frac{1}{2\pi} \int_{\mathbb{S}^2} \chi_{a,b} \left(x\cdot \omega + t\right) G(x\cdot \omega +t, \omega) d\omega
 \]
 Here $\chi_{a,b}$ is the characteristic function of $(a,b)$. We may also write $u(x,t)$ in term of $\mathbf{T}$
 \begin{align*}
  &u(x,t) = \frac{1}{2\pi} \mathbf{T} \tilde{G}^{(t)},& &\tilde{G}^{(t)} (s, \omega)= (\mathbf{P}_{a,b} G) (s+t,\omega) = \chi_{a,b} \left(s + t\right)  G(s +t, \omega).&
 \end{align*}
 Since $\tilde{G}^{(t)}$ is supported in the set $(-t + a, -t + b)\times \mathbb{S}^2$, we may apply the $L^6$ decay estimate \eqref{left point interpolation} and obtain
 \[
  \|u(\cdot,t)\|_{L^6(\{x: |x|>R\})} \lesssim \left(\frac{b-a}{R}\right)^{1/3} \|\tilde{G}^{(t)}\|_{L^2(\Rm \times \mathbb{S}^2)} \leq \left(\frac{b-a}{R}\right)^{1/3} \|G\|_{L^2(\Rm \times \mathbb{S}^2)}.
 \]
 In addition, if $t$ tends to infinity, we may apply the global estimate \eqref{left point interpolation 2} to obtain
 \[
  \|u(\cdot,t)\|_{L^6(\Rm^3)} \lesssim |t|^{-1/3} \|G\|_{L^2(\Rm \times \mathbb{S}^2)}.
 \]
 Therefore the map defined above is a bounded linear operator from $L^2(\Rm \times \mathbb{S}^2)$ to the space $L_0^\infty (\Rm_t; L^6(\{x: |x|>R\}))$ whose operator norm is dominated by $(b-a)^{1/3} R^{-1/3}$. Here the space $L_0^\infty$ is a subspace of the regular $L^\infty$ space defined in the following way: We first consider the space of all simple functions so that each simple function is zero except in a set with a finite measure then define the completion of this space in the $L^\infty$ space to be $L^\infty_0$. Next we consider 
\begin{align*}
  &G \in \dot{H}^{s_0(p)-1} (\Rm \times \mathbb{S}^2),& &s_0(p) = \frac{5p-9}{(9-p)(p-1)} \in (1/2, s_p).&
\end{align*}
By the Strichartz estimates and Proposition \ref{prop isometry}, we have 
\begin{align*}
 \|u(x,t)\|_{L^\frac{2(9-p)(p-1)}{5p-9} L^\frac{2(9-p)(p-1)}{p(5-p)} (\Rm \times \Rm^3)} &\lesssim \|\mathbf{T}_-^{-1} \mathbf{P}_{a,b} G\|_{\dot{H}^{s_0(p)}\times \dot{H}^{s_0(p)-1}(\Rm^3)}\\
 & \lesssim \|\mathbf{P}_{a,b} G\|_{\dot{H}^{s_0(p)-1} (\Rm \times \mathbb{S}^2)} \\
 & \lesssim \|G\|_{\dot{H}^{s_0(p)-1}(\Rm \times \mathbb{S}^2)}.
\end{align*}
Thus the map defined above is also a bounded operator from $\dot{H}^{s_0(p)-1}(\Rm \times \mathbb{S}^2)$ to the space $L^\frac{2(9-p)(p-1)}{5p-9}(\Rm_t; L^\frac{2(9-p)(p-1)}{p(5-p)}(\{x: |x|>R\}))$ whose operator norm is smaller or equal to a constant solely determined by $p$. We then apply the complex interpolation method and conclude that 
\[
 \|u(x,t)\|_{L^\frac{2p}{1+s_p} L^\frac{2p}{2-s_p}(\Rm \times \{x: |x|>R\})} \lesssim \left(\frac{b-a}{R}\right)^{\frac{p-3}{2p}}\|G\|_{\dot{H}^{s_p-1}(\Rm \times \mathbb{S}^2)}.
\]
Please see \cite{interpolationspace} for the details of interpolation method and spaces, for instance. Finally we observe that $\mathbf{P}_{a,b} G = G$ if $G \in \dot{H}_0^{s_p-1}(I \times \mathbb{S}^2)$ and thus finish the proof. 
\end{proof}
\section{Non-radiative solutions to non-linear equations} \label{sec: non-linear wave equation}
In this section we assume $p \in (3,5)$ and consider non-linear wave equations
\[
 \left\{\begin{array}{ll} \partial_t^2 u - \Delta u = F(x,t,u); & (x,t) \in \Rm^3 \times \Rm; \\ (u,u_t)|_{t=0} = (u_0,u_1)\in \mathcal{H}^{s_p}. \end{array}\right.
\]
Here $s_p = 3/2-2/(p-1)$, the non-linear term $F(x,t,u)$ satisfies 
\begin{align*}
 &|F(x,t,u)| \lesssim |u|^p;& &|F(x,t,u_1)-F(x,t,u_2)| \lesssim (|u_1|^{p-1} +|u_2|^{p-1}) |u_1-u_2|.&
\end{align*}
We recall the definition of $Y(I)$ and $Z(I)$ norms. A basic calculation shows that 
\begin{align*}
 \|F(x,t,u)\|_{Z(I)} &\lesssim \|u\|_{Y(I)}^p;\\
  \|F(x,t,u_1)-F(x,t,u_2)\|_{Z(I)} & \lesssim (\|u_1\|_{Y(I)}^{p-1} + \|u_2\|_{Y(I)}^{p-1})\|u_1-u_2\|_{Y(I)}.
\end{align*}
Here $I$ is an arbitrary time interval. Combining these inequalities with the Strichartz estimates, we may apply a fixed-point argument and verify the well-posedness of this equation. More precisely, there exists a maximal interval of existence $I$ and a unique solution defined for $t \in I$ so that 
\begin{itemize}
 \item $(u,u_t) \in C(\Rm_t; \mathcal{H}^{s_p})$; 
 \item $u \in Y(J)$ for any bounded closed interval $J \subset I$. 
\end{itemize}
Next we give the definition of non-radiative solutions, in the energy sub-critical setting. 
 \begin{definition}
 We say a solution to be an $R$-weakly non-radiative solution, if there exist a sequence of time $t_k$ and data $(v_{0,k}, v_{1,k})$ with the limits 
 \begin{align*}
  &\lim_{k\rightarrow \pm \infty} t_k = \pm \infty;& &\lim_{k\rightarrow \pm \infty} \|(v_{0,k}, v_{1,k})\|_{\dot{H}^{s_p}\times \dot{H}^{s_p-1}(\Rm^3)} = 0;&
 \end{align*}
 so that $(u(x, t_k), u_t(x,t_k))= (v_{0,k}(x), v_{1,k}(x))$ holds in the exterior region $\{x: |x|>R+|t_k|\}$. 
\end{definition}
\begin{remark}
 A soliton-like solution, i.e. a solution to (CP1) so that its orbit 
 \[
  \{(u(\cdot,t), u_t(\cdot,t)): t \in \Rm\}
 \]
 is pre-compact in $\mathcal{H}^{s_p}$, must be non-radiative. This kind of solutions are frequently used in the compact-rigidity argument. Please see, for example, Dodson-Lawrie-Mendelson-Murphy \cite{nonradial3p5}. 
\end{remark}
\begin{lemma} \label{non-radiative to G}
 Assume that $u(x,t)$ is an $R$-weakly non-radiative solution and $v$ is a free wave, so that $u$ scatters to $v$ in the positive time direction:
 \[
  \lim_{t\rightarrow +\infty} \|(u(\cdot,t), u_t(\cdot, t)) - (v(\cdot,t), v_t(\cdot,t))\|_{\dot{H}^{s_p}\times \dot{H}^{s_p-1}(\Rm^3)} = 0. 
 \]  
 Let $G^+(s,\omega)$ be the radiation profile of $v$ in the positive time direction. Then we must have $G^+(s,\omega) = 0$ in the region $(R,+\infty) \times \mathbb{S}^2$. 
\end{lemma}
\begin{remark}
 A similar lemma holds in the negative time direction since the wave equation is time-reversible. 
\end{remark}
\begin{proof}
 It suffices to show that 
 \[
  \int_{\Rm \times \mathbb{S}^2} G^+(s,\omega) \rho (s,\omega) ds d\omega = 0
 \]
 for any test function $\rho(s,\omega) \in C_0^\infty ((R, +\infty)\times \mathbb{S}^2)$. We first combine our radiation fields and the scattering assumption to obtain 
 \begin{equation} \label{u close Gplus}
   \lim_{t\rightarrow +\infty} \left\|u_t(x,t) - \frac{1}{|x|} G^+(|x|-t, x/|x|)\right\|_{\dot{H}^{s_p-1}(\Rm^3)} = 0
 \end{equation}
 Given any $\rho(s, \omega) \in C_0^\infty((R,+\infty)\times \mathbb{S}^2)$, a basic calculation shows that the functions
 \[
   \rho^{(t)} (x) = \frac{1}{|x|} \rho(|x|-t, x/|x|) 
 \]
 is uniformly bounded in the space $H^1(\Rm^3)$ for all time $t \gg 1$. Thus they are also uniformly bounded in the space $\dot{H}^{1-s_p}(\Rm^3)$. We may combine this with \eqref{u close Gplus} and obtain
 \begin{align*}
  \lim_{t\rightarrow +\infty} \int_{\Rm^3} u_t(x,t) \rho^{(t)}(x) dx & = \lim_{t\rightarrow +\infty} \int_{\Rm^3}  \frac{1}{|x|} G^+(|x|-t, x/|x|) \rho^{(t)}(x) dx\\
  & = \int_{\Rm \times \mathbb{S}^2} G^+(s,\omega) \rho(s, \omega) ds d\omega.
 \end{align*}
 Finally we observe that the non-radiative assumption on $u$ implies the limit of the left hand side must be zero and finish the proof.
\end{proof}

\paragraph{Power-like decay by recursion formula} Before we give the main result in this section, we first recall a technical lemma, which will be used in the proof of our main result. This kind of results have been known and used for a long time. Please see, for example, Li-Shen-Wang-Wei \cite{nonradialCE} for a  proof.
\begin{lemma} \label{recursion lemma}
Assume that $l>1$ and $\alpha>0$ are constants. Let $S: [R,+\infty) \rightarrow [0,+\infty)$ be a function satisfying 
\begin{itemize}
 \item $S(r)\rightarrow 0$ as $r\rightarrow +\infty$;
 \item The recursion formula $S(r_2) \lesssim (r_1/r_2)^\alpha + S^l (r_1)$ holds when $r_2 \gg r_1 \gg R$.
\end{itemize}
Then given any constant $\beta \in (0, (1-1/l)\alpha)$, the decay estimate $S(r) \leq r^{-\beta}$ holds as long as $r>R_0$ is sufficiently large. 
\end{lemma}

\noindent Now we introduce our main result

\begin{proposition} \label{main 2}
 Assume that $u$ is an $R$-weakly non-radiative solution as defined above. Let $G$ be the radiation field of the initial data $(u_0,u_1)$ in the negative time direction. Then for sufficiently large $r > R_0(u)$, we have 
 \begin{align*}
  \|\mathbf{P}_r^\pm G\|_{\dot{H}^{s_p-1}(\Rm \times \mathbb{S}^2)} \lesssim r^{-\frac{p-3}{2}}.
 \end{align*} 
\end{proposition}
\begin{proof}
Given $r \geq R$, we consider the modified non-linear wave equation
\[
 \left\{\begin{array}{ll} \partial_t^2 v - \Delta v = \chi_r(x,t) F(x,t,u), & (x,t)\in \Rm^3 \times \Rm; \\
 (u,u_t)|_{t=0} = (u_0,u_1). & \end{array}\right.
\]
Here $\chi_r(x,t)$ is the characteristic function of the exterior region $\Omega_r \doteq \{(x,t): |x|>|t|+r\}$. First of all, the following inequalities hold by our assumption on the nonlinear term $F$. 
\begin{align*}
 \left\|\chi_r F(x,t,v)\right\|_{Z (\Rm)} & \lesssim \|\chi_r v\|_{Y(\Rm)}^p; \\
 \left\|\chi_r F(x,t,v_1) - \chi_r F(x,t,v_2)\right\|_{Z(\Rm)} & \lesssim (\|\chi_r v_1\|_{Y(\Rm)}^{p-1} + \|\chi_r v_2\|_{Y(r)}^{p-1})\|\chi_r v_1 - \chi_r v_2\|_{Y(\Rm)}.
\end{align*}
We may apply a fixed-point argument (see, for instance, \cite{locad1, ls, Pecher}) and conclude that as long as $S(r) \doteq \|\chi_r \mathbf{S}_L (u_0,u_1)\|_{Y(\Rm)}$ is sufficiently small, which always holds if $r$ is sufficiently large, then the modified wave equation above has a unique solution $v \in C(\Rm_t; \dot{H}^{s_p}\times \dot{H}^{s_p-1}) \cap Y(\Rm)$, with 
\begin{align}
 \|\chi_r  v\|_{Y(\Rm)} & \leq 2 S(r); \\
 \sup_{t} \|(v(\cdot,t),v_t(\cdot,t))\|_{\mathcal{H}^{s_p}} & \leq \|(u_0,u_1)\|_{\mathcal{H}^{s_p}} + C(p) S^p (r). \label{sup bound v}
\end{align}
In addition, there exists two free waves $v_1(x,t)$ and $v_2(x,t)$ so that $v$ scatters to $v_1$, $v_2$ as $t\rightarrow - \infty$ and $t\rightarrow +\infty$, respectively. More precisely we have
\begin{align*}
 \lim_{t\rightarrow -\infty} \|(v(\cdot,t),v_t(\cdot,t))-v_1(\cdot,t),\partial_t v_1(\cdot,t))\|_{\mathcal{H}^{s_p}} & = 0;\\
  \lim_{t\rightarrow +\infty} \|(v(\cdot,t),v_t(\cdot,t))-v_2(\cdot,t),\partial_t v_2(\cdot,t))\|_{\mathcal{H}^{s_p}} & = 0.
\end{align*}
For convenience we use the notations $G_1, G_2$ for the radiation profiles of $v_1, v_2$ in the negative time direction, respectively. By finite speed of propagation, we have $u(x,t) = v(x,t)$ in the region $\Omega_r$. Thus $v(x,t)$ is an $r$-weakly non-radiative solution. According to Lemma \ref{non-radiative to G}, we have 
\begin{align}
 &G_1(s,\omega) = 0, \quad s>r;& &G_2(s,\omega) = 0, \quad s<-r.& \label{support of G12}
\end{align}
Next we observe ($u_L = \mathbf{S}_L (u_0,u_1)$)
\begin{align*}
 \|G_1 - G\|_{\dot{H}^{s_p-1}(\Rm \times \mathbb{S}^2)} &= \frac{1}{\sqrt{2}} \lim_{t\rightarrow -\infty} \|(v_1(\cdot,t), \partial_t v_1(\cdot,t)) - (u_L (\cdot,t), \partial_t u_L(\cdot,t))\|_{\mathcal{H}^{s_p}}\\
 & =  \frac{1}{\sqrt{2}} \lim_{t\rightarrow -\infty} \|(v(\cdot,t), \partial_t v(\cdot,t)) - (u_L (\cdot,t), \partial_t u_L(\cdot,t))\|_{\mathcal{H}^{s_p}}\\
 & \lesssim_p \|\chi_r F(x,t,v)\|_{Z((-\infty, 0])}\\
 & \lesssim_p \|\chi_r v\|_{Y(\Rm)}^p\\
 & \lesssim_p S^p(r).
 \end{align*}
 Similarly we have 
 \[
  \|G_2 - G\|_{\dot{H}^{s_p-1}(\Rm \times \mathbb{S}^2)} \lesssim_p S^p (r).
 \]
 We then write $G$ in a sum
 \[
  G = \mathbf{P}_{-r,r} G + \mathbf{P}_r^+ G + \mathbf{P}_r^- G. 
 \]
Here we have
\[
 \|\mathbf{P}_{-r,r} G\|_{\dot{H}_0^{s_p-1}((-r, r)\times \mathbb{S}^2)}  \lesssim_p  \|G\|_{\dot{H}^{s_p-1}}.
\]
In addition, we recall \eqref{support of G12} and obtain
\begin{align*}
 \|\mathbf{P}_r^+ G\|_{\dot{H}^{s_p-1}(\Rm \times \mathbb{S}^2)} & =   \|\mathbf{P}_r^+ (G-G_1)\|_{\dot{H}^{s_p-1}(\Rm \times \mathbb{S}^2)}  \\
 & \lesssim \|G-G_1\|_{\dot{H}^{s_p-1}(\Rm \times \mathbb{S}^2)} \\
 & \lesssim S^{p} (r).
\end{align*}
Similarly we have 
\[
 \|\mathbf{P}_r^- G\|_{\dot{H}^{s_p-1}(\Rm \times \mathbb{S}^2)} \lesssim_p S^p (r). 
\]
Thus we may apply the regular Strichartz estimates, as well as the decay estimates given in Lemma \ref{Y estimate} to obtain ($r_1>r$)
\begin{align*}
  \|\chi_{r_1} u_L\|_{Y(\Rm)} & \leq \|u_L\|_{L^\frac{2p}{1+s_p}(\Rm_t; L^\frac{2p}{2-s_p}(\{x: |x|>r_1\}))}\\
 & \lesssim_p \left(r/r_1\right)^\frac{p-3}{2p} \|\mathbf{P}_{-r,r} G\|_{\dot{H}_0^{s_p-1}} + 
  \|\mathbf{P}_r^+ G + \mathbf{P}_r^- G\|_{\dot{H}^{s_p-1}(\Rm \times \mathbb{S}^2)} \\
  & \lesssim_p S^p(r) + (r/r_1)^\frac{p-3}{2p} \|G\|_{\dot{H}^{s_p-1}}.
\end{align*}
Namely we have 
\[
 S(r_1) \lesssim_p S^p (r) + (r/r_1)^\frac{p-3}{2p} \|G\|_{\dot{H}^{s_p-1}}, \qquad r_1>r\gg R; 
\]
with $S(r) \rightarrow 0$ as $r \rightarrow +\infty$. It immediately follows that given any constant $\alpha \in (0, \frac{(p-3)(p-1)}{2p^2})$, there exists $R_1 = R_1(\alpha,u) \gg R$, so that
\[
 S(r) \leq r^{-\alpha}, \qquad r\geq  R_1.   
\]
Thus we also have 
\[
 \|\mathbf{P}_r^\pm G\|_{\dot{H}^{s_p-1}} \lesssim_p r^{-p \alpha}, \qquad r\geq R_1. 
\]
We may write 
\begin{equation} \label{decomposition of G}
 G =  G_0 + \sum_{k=1}^\infty (G_k^+ + G_k^-).
\end{equation} 
Here 
\begin{align*}
 &G_0 = \mathbf{P}_{-R_1,R_1} G;& &G_k^+ = \mathbf{P}_{2^{k-1} R_1, 2^k R_1} G;& &G_k^- = \mathbf{P}_{-2^{k} R_1, -2^{k-1} R_1} G.&
\end{align*}
We have $\|G_k^\pm\|_{\dot{H}^{s_p-1}} \lesssim \|\mathbf{P}_{2^{k-1} R_1}^\pm G\|_{\dot{H}^{s_p-1}} \lesssim 2^{-p\alpha k}$. This immediately gives a stronger upper bound of $S(r)$ by Lemma \ref{Y estimate}.
\begin{align*}
 S(r) &\lesssim_p \left(R_1/r\right)^\frac{p-3}{2p}\|G_0\|_{\dot{H}^{s_p-1}} +  \sum_{k=1}^{+\infty} \sum_\pm \left(2^{k}R_1/r\right)^\frac{p-3}{2p}\|G_k^\pm\|_{\dot{H}^{s_p-1}} \\
 & \lesssim r^{-\frac{p-3}{2p}} + \sum_{k=1}^{+\infty} 2^{\left(\frac{p-3}{2p} - p\alpha\right)k} r^{-\frac{p-3}{2p}}.
\end{align*}
We may choose $\alpha = \alpha(p)  \in (\frac{p-3}{2p^2}, \frac{(p-3)(p-1)}{2p^2})$ and obtain
\[
 S(r) \lesssim r^{-\frac{p-3}{2p}}.
\]
Finally we recall the estimates $\|\mathbf{P}_r^\pm G\|_{\dot{H}^{s_p-1}(\Rm \times \mathbb{S}^2)} \lesssim_p S^p(r)$ and conclude
\[
 \|\mathbf{P}_r^\pm G\|_{\dot{H}^{s_p-1}(\Rm \times \mathbb{S}^2)} \lesssim r^{-\frac{p-3}{2}}, \qquad r \gg R.
\]
\end{proof}

\begin{corollary} \label{lower s}
 If $u$ is a weakly non-radiative solution as defined above, then its initial data $(u_0,u_1)$ satisfy
 \[
  (u_0,u_1) \in \dot{H}^{\beta} (\Rm^3) \times \dot{H}^{\beta-1} (\Rm^3), \qquad \forall \beta \in (1/2, s_p). 
 \]
\end{corollary}
\begin{proof}
 First of all, we recall that given $a<b$, then we have 
 \[
  \|\mathbf{P}_{a,b} f(s)\|_{\dot{H}^{\beta-1}(\Rm)} \lesssim \|f\|_{\dot{H}^{\beta-1}(\Rm)}
 \]
 Here the implicit constant is independent to $a,b$. In addition, we may utilize the Sobolev embedding and obtain 
 \[
  \|\mathbf{P}_{a,b} f(s)\|_{\dot{H}^{\beta-1}(\Rm)} \lesssim \|\mathbf{P}_{a,b} f(s)\|_{L^{\frac{1}{3/2-\beta}}(\Rm)} \lesssim (b-a)^{1-\beta} \|f\|_{L^2(\Rm)}.
 \]
 An interpolation shows that 
 \[
  \|\mathbf{P}_{a,b} f(s)\|_{\dot{H}^{\beta-1}(\Rm)} \lesssim (b-a)^{s_p-\beta} \|f\|_{\dot{H}^{s_p-1}(\Rm)}. 
 \]
 Thus we have 
 \begin{equation} \label{Sobolev embedding}
  \|\mathbf{P}_{a,b} f(s,\omega)\|_{\dot{H}^{\beta-1}(\Rm\times \mathbb{S}^2)} \lesssim (b-a)^{s_p-\beta} \|f\|_{\dot{H}^{s_p-1}(\Rm\times \mathbb{S}^2)}.
 \end{equation} 
 Let $G$ be the radiation profile of $(u_0,u_1)$ in the negative time direction. We recall Proposition \ref{main 2} and utilize the decomposition of $G$ as given in \eqref{decomposition of G}, where $R_1$ is chosen so that $\|\mathbf{P}_r^\pm\|_{\dot{H}^{s_p-1}} \lesssim r^{-\frac{p-3}{2}}$ if $r\geq R_1$. We may apply the local Sobolev embedding \eqref{Sobolev embedding} and obtain 
 \[
  \|G_0\|_{\dot{H}^{\beta-1}} = \|\mathbf{P}_{-R_1,R_1} G\|_{\dot{H}^{\beta-1}} \lesssim R_1^{s_p - \beta} \|G\|_{\dot{H}^{s_p-1}} < +\infty.
\]
We also have
 \[
  \|G_k^+\|_{\dot{H}^{\beta-1}} = \left\|\mathbf{P}_{2^{k-1} R_1, 2^k R_1} \mathbf{P}_{2^{k-1} R_1}^+ G\right\|_{\dot{H}^{\beta-1}} \lesssim (2^k R_1)^{s_p-\beta} \|\mathbf{P}_{2^{k-1} R_1}^+ G\|_{\dot{H}^{s_p-1}} \lesssim (2^k R_1)^{s_p-\beta - \frac{p-3}{2}}. 
 \]
 Similarly we have 
 \[
   \|G_k^-\|_{\dot{H}^{\beta-1}} \lesssim (2^k R_1)^{s_p-\beta - \frac{p-3}{2}}. 
 \]
 The series 
 \[
  \sum_{k=1}^\infty (2^k R_1)^{s_p-\beta - \frac{p-3}{2}}
 \]
 converges by the fact 
 \[
  s_p -\beta - \frac{p-3}{2} < s_p - \frac{1}{2} - \frac{p-3}{2} = - \frac{(p-3)^2}{2(p-1)} < 0.
 \]
 As a result, we have $G \in \dot{H}^{\beta-1}(\Rm)$. Finally we apply Proposition \ref{prop isometry} to conclude that $(u_0,u_1) \in \mathcal{H}^\beta (\Rm^3)$. 
\end{proof}
\begin{remark}
 A similar result holds in the energy critical case as well. More precisely, let $u$ be an $R$-weakly non-radiative solution to the non-linear equation 
 \[
  \left\{\begin{array}{ll} \partial_t^2 u - \Delta u = F(x,t,u), & (x,t) \in \Rm^3 \times \Rm; \\
  u(\cdot,0) = u_0 \in \dot{H}^1 (\Rm^3); & \\
  u_t(\cdot,0) = u_1 \in L^2 (\Rm^3); \end{array}\right.
 \]
 whose non-linear term $F(x,t,u)$ satisfies 
 \begin{align*}
  &|F(x,t,u)| \lesssim |u|^5;& &|F(x,t,u_1)-F(x,t,u_2)| \lesssim (|u_1|^4 + |u_2|^4)|u_1-u_2|.&
 \end{align*}
 Then we have $(u_0,u_1) \in \dot{H}^{\beta} \times \dot{H}^{\beta-1}(\Rm^3)$ for any $\beta \in (1/2, 1)$. Our previous work \cite{nonradialCE} shows that the radiation profile $G$ of initial data satisfies $\|\mathbf{P}_r^\pm G\|_{L^2(\Rm \times \mathbb{S}^2)} \lesssim r^{-\alpha}$ for any $\alpha \in (0,1)$ and sufficiently large $r\geq R(u,\alpha)$. We then follow exactly the same argument as Corollary \ref{lower s}. 
\end{remark}
\section*{Acknowledgement}
The authors are financially supported by National Natural Science Foundation of China Project 12071339.

\end{document}